\documentclass[11pt]{article}
\usepackage{amsmath,amssymb}
\usepackage{tikz-cd}
\newtheorem{thm}{Theorem}[section]
 \newtheorem{cor}[thm]{Corollary}

 \newtheorem{defn}[thm]{Definition}
 \newtheorem{rem}[thm]{Remark}
 \newtheorem{ex}{Example}
  \newtheorem{pr}{Problem}
 \numberwithin{equation}{section}
\usepackage{geometry}
\begin{document}

% ------------------------------------------------------------------------

%-------------------------------------------------------------------------
% editorial commands: to be inserted by the editorial office
%
%\firstpage{1} \volume{228} \Copyrightyear{2004} \DOI{003-0001}
%\seriesextraline{This is the Concrete Title of this Book\br H.E. R and S.T.C. W, Eds.}
%\firstpage{1}

\begin{center}
\textbf{\Large Riemann integrability under weaker forms of continuity in infinite dimensional space}
\end{center} 

\vspace{1cm}

%----------Author 1
\begin{center}
\textbf{M.A. Sofi}\\
\textit{Department of Mathematics\\
University of Kashmir, Srinagar-190006, India}\\
\end{center}
 \begin{center}
 email: aminsofi@gmail.com
 \end{center}


\date{January 1, 2004}

%%% ----------------------------------------------------------------------

\begin{abstract}
In classical analysis, the relationship between continuity and Riemann integrability of a function is an intimate one: a continuous function $f:[a,b]\rightarrow \mathbb{R}$ is always Riemann integrable whereas a Riemann integrable function is continuous almost everywhere (a.e). In the setting of functions taking values in infinite dimensional abstract spaces that include quasi Banach spaces, one encounters certain curious phenomena involving the breakdown of the above stated phenomena, besides the failure of the fundamental theorem of calculus and the non-existence of primitives for continuous functions! While these properties can surely be salvaged within the class of Banach spaces, it turns out that certain important properties involving vector integration that include  Riemann integration no longer hold in an infinite dimensional setting. This will be seen to be the case, for example, in situations when it is required to integrate functions which are continuous with respect to certain well known linear topologies on $X$ (resp. on $X^{*}$) weaker than the norm topology. As we shall see in Section $3\text{‬(a)},$ such a requirement imposes rather severe restrictions on the space in question. The present paper is devoted to a discussion of these issues which will be examined in the setting of Banach and Fr\'{e}chet spaces on the one hand and of quasi Banach spaces on the other.
\end{abstract}
 \textbf{Subject Class [2010]:} Primary 46G10; Secondary  28B05
\section{Riemann integrability in Banach spaces}\label{s1}
Let us recall that a (bounded) function $f:[a,b]\rightarrow X$ defined on a closed and bounded interval and taking values in a Banach space $X$ is said to be Riemann integrable if the following condition $(*)$ holds:

$(*)$ $\exists\, x\in X$ such that $\forall\, \epsilon>0\, \exists\,\delta=\delta(\epsilon)>0$ such that for each (tagged) partition $P= \{s_i,[t_{(i-1)},t_i ],1\leq i\leq j\}$ of $[a, b]$ where $a=t_0<t_1<\cdots<t_j=b$ and $s_i\in[t_{(i-1)},t_i ],1\leq i\leq j,$ with
\begin{equation}\label{1.1}
\|P\|=\underset{1\leq j\leq n}{\max}(t_{j}-t_{j-1})<\delta,
\end{equation}
we have
\begin{equation*}
\|S(f,P,\xi)-x\|\leq \epsilon,
\end{equation*}
where $\xi=\{s_i;1\leq i\leq j\}$ and $ S(f,P,\xi)$ is the Riemann sum of $f$ corresponding to the partition  $P$:
\begin{equation*}
S(f,P,\xi)=\sum\limits_{i=1}^{j}f(s_i)(t_i-t_{i-1}).
\end{equation*}
The (unique) vector $x,$ to be denoted by 
$$\int\limits_{a}^{b}f(t)dt  $$
shall be called the Riemann-integral of $f$ over $[a, b].$ 

\noindent(The partition $P$ satisfying \eqref{1.1} shall be referred to as a $\delta$-fine partition). 

The same definition holds for $X$-valued functions where $X$ is a quasi Banach space in which case the norm shall be replaced by the corresponding quasi norm on $X.$ A slight modification of this definition shall provide an appropriate analogue of Riemann integrability in Fr\'{e}chet spaces which we shall pursue in Section \ref{s3}(c).

In this section, we shall study the abstract analogues of the following results which are well known from classical analysis: 
\begin{enumerate}
\item[a.] Given a continuous function $f:[a,b]\rightarrow\mathbb{R},$ it holds that 
\begin{enumerate}
\item[(i)] $f$ is Riemann integrable.
\item[(ii)] $f$ has a primitive: $\exists$ $F:[a,b]\rightarrow \mathbb{R}$ such that $F$ is differentiable on $[a,b]$ and
$$ F^{\prime}(t)=f(t)\quad\text{on}\quad [a,b]. $$
\end{enumerate}
\item[b.] (Fundamental Theorem of Calculus): Given that $f:[a,b]\rightarrow\mathbb{R}$ is differentiable and that $f^{\prime}$ is Riemann integrable, we have
$$ f(x)=\int\limits_{a}^{x}f^{\prime}(t)dt,\qquad\forall x\in [a,b]. $$
\item[c.] Let $f:[a,b]\rightarrow\mathbb{R}$ be a bounded function. Then $f$ is Riemann integrable, if and only if it is continuous almost everywhere.

In the case of Banach space valued functions $f:[a,b]\rightarrow X,$ it turns out that the same method of proof as employed in the case of real valued functions can be slightly modified to yield the following Banach space analogues of these results: 
\item[d.] Let $X$ be a Banach space and  $f:[a,b]\rightarrow X$ a continuous function. Then 
\begin{enumerate}
\item[(i)] $f$ is Riemann integrable.
\item[(ii)] $f$ has a primitive: $\exists$ $F:[a,b]\rightarrow X$ such that $F$ is differentiable on $[a,b]$ and
$$ F^{\prime}(t)=f(t)\quad\text{on}\quad [a,b]. $$
In fact, we can choose $F$ to be given by
$$ F(x)=\int\limits_{a}^{x}f(t)dt,\qquad x\in [a,b]. $$
\end{enumerate}
\item[e.] (Fundamental Theorem of Calculus): If $f$ is differentiable on $[a, b],$ then $f^{\prime}$ is always Henstock integrable and
$$ f(x)=\int\limits_{a}^{x}f^{\prime}(t)dt,\qquad\forall x\in [a,b]. $$
Here, Henstock integrability is meant exactly in the sense of Riemann integrability -as defined above- except that $\delta$ appearing in the definition of Riemann integrability has to be replaced by a positive function on $\mathbb{R}$ (called a \textit{gauge}), the rest remaining the same, with \eqref{1.1} now being replaced by
$$ t_{i}-t_{i-1}\leq \delta(s_i),\quad i=1,2,\cdots,n. $$
\begin{rem}
(i) Choosing the gauge $\delta$ to be a constant function, it is clear from the definition that Riemann integration theory is subsumed within the more general framework of Henstock integration theory and, therefore, (e) yields (FTC) for Riemann integrable functions, provided it is assumed that $f^{\prime}$ is Riemann integrable.\\
(ii) It is remarkable that the simple idea of replacing $\delta$ by a positive function leads to a powerful integration theory in which all the convergence theorems of Lebesgue integration theory hold and more importantly, the fundamental theorem of calculus (FTC) holds in its utmost generality without assuming the integrability of the derivative, as seen in (e) above. In this sense, Henstock integration theory scores over Lebesgue integration in view of the simplicity of the (FTC) in this new setting.
\end{rem}
\item[f.] If $f$ is continuous a.e., then it is Riemann integrable. However, as we shall see below, converse is not true!
\end{enumerate}
The assertion in (f) above motivates the following definition.
\begin{defn}[\cite{7}]\label{d1.1}
A Banach space $X$ is said to have Riemann-Lebesgue (RL)- property  if each (bounded) Riemann integrable function $f:[a,b]\rightarrow X$ is continuous a.e on $[a, b].$
\end{defn}

\begin{ex}\label{e1}
\begin{enumerate}
\item[(a)] (Lebesgue):  $\mathbb{R}$ has (RL)-property.\\
Consequence: Finite dimensional Banach spaces have the (RL)-property.
\item[(b)] (G.C.da Rocha \cite{17}): $\ell_1$ has (RL)-property.
\item[(c)] (G.C.da Rocha \cite{17}):Tsirelson space has (RL)-property.
\end{enumerate}
\end{ex}
Before proceeding further, let us pause to see how the proof of (b) works. We shall slightly modify the proof as given in \cite{17} (see also \cite{7}) to prove the more general statement that $\ell_1(X)$- the space of $X$-valued absolutely summable sequences - has (RL)-property whenever $X$ has it. Thus, let $f:[0,1]\rightarrow \ell_1 (X)$ be Riemann integrable and assume, on the contrary, that $f$ is discontinuous on a set $H$ of positive measure. For a given $t\in[a,b],$ denote by $\omega(f,t)$ the oscillation of $f$ at $t$: 
$$ \omega(f,t)=\underset{\epsilon\rightarrow 0}{\lim}\,\omega(f,[t-\epsilon,t+\epsilon]) $$
where $\omega(f,[a,b])$ is defined by
$$ \omega(f,[0,1])=\sup\{\|f(x)-f(y)\|: x,y\in [0,1]\}. $$
Thus, $f$ is continuous at $t$ if and only if $\omega(f,t)=0$ and so, there exist $\alpha>0,$ $\beta>0$ such that $\mu(H)=\alpha$ and $H=\{t\in[0,1];\omega(f,t)\geq \beta \},$ where $\mu$ denotes the Lebesgue measure. To get a contradiction, it suffices to show that for each $\delta>0,$ there exist tagged partition $(P_1,\xi_1)$ and $(P_2,\xi_2)$ of $[0,1]$ such that $\|P_1\|<\delta,$ $\|P_2\|<\delta$ with 
$$\|S(f,P_1,\xi_1)-S(f,P_1,\xi_1)\|>\alpha\beta/4.\eqno{(**)}$$
For each $n\geq 1,$ denote by $e_{n}^{*}$ the $n$th co-ordinate functional on $\ell_1(X):e_{n}^{*} (\{x_i \})=x_n$ and let $G_n$ be the set of all discontinuities of $e_{n}^{*}\circ f.$  Assuming that $\mu(G_n)\neq 0$ for some $n\geq 1$ and combining this with the (RL)-property of $X$ gives that $e_{n}^{*}\circ f$  and hence $f$ is not Riemann integrable, contradicting the given hypothesis. This gives $\mu(G_n)=0$ for each $n\geq 1$ and, therefore, the set $G=\bigcup_{n=1}^{\infty}G_n$  has measure zero at each point of which $e_{n}^{*}\circ f$ is continuous for each $n\geq 1.$
	
Fix  $\delta>0$ and choose $N$ such that $1/N<\delta.$ Let $D=\{k/N;k=0,…N)\}$ and denote by $P_1$ and $P_2$ the partitions of $[0,1]$ determined by the points of $D.$ To specify the tags $\xi_1$ and $\xi_2$ of these partitions, we let $P$ be the collection of all intervals $[a_i,b_i ]$ from (either of) these partitions such that $\mu(H\bigcap(a_i,b_i))>0.$ If $p$ is the number of elements in $P,$ it follows that $p\geq N\alpha.$ For each $i=1, . . . p,$ fix $u_i\in (H\setminus G)\bigcap[a_i,b_i],$ arbitrarily. (This is possible because $G$ has measure zero).  We shall construct sets $\{v_i;1\leq i\leq p \}\subset[a_i,b_i ]$ and $\{0=n_0\leq n_1\leq \cdots\leq n_p\}\subset\mathbb{N}$ satisfying the following properties- with $z^{(i)}$ defined to be the element of $\ell_1 (X)$ given by: $z^{(i)}=f(u_i )-f(v_i ),\,1\leq i\leq p.$
\begin{enumerate}
\item[(a)] $\left\|z^{(i)}\right\|\geq \beta/2$
\item[(b)] $\sum_{j=n_{i}}^{\infty}\left\|z_{j}^{(i)}\right\|<\epsilon 2^{-i}$
\item[(c)] $\sum_{j=1}^{n_{i-1}}\left\|z_{j}^{(i)}\right\|<\epsilon 2^{-i}.$
\end{enumerate}
Claim: $\left\|\sum_{j=1}^{p}z^{(i)}\right\|\geq \left(\frac{p\beta}{2}-4\epsilon\right).$

\noindent For each $1\leq i\leq p,$ consider $y^{(i)}\in\ell_i(X)$ given by
\begin{align}
y_{j}^{(i)}=
\begin{cases}
z_{j}^{(i)}, \,&n_{i-1}<j\leq n_{i}\\
0,\, &\text{otherwise}
\end{cases}
\end{align}
We note that for each $1\leq i\leq p,$ we have
\begin{align*}
\|z^{(i)}-y^{(i)}\|=\sum\limits_{j=1}^{n_{i-1}}\|z_j^{(i)}\|+\sum\limits_{j=n_{(i)}}^{\infty}\|z_j^{(i)}\|<2\epsilon 2^{-i}
\end{align*}
and therefore
\begin{align*}
\|y^{(i)}\|\geq \|z^{(i)}\|-\|z^{(i)}-y^{(i)}\|>\beta/2-2\epsilon 2^{-i}.
\end{align*}
Combining these estimates yields
\begin{align*}
\left\|\sum\limits_{i=1}^{p}z^{(i)}\right\|&\geq \left\|\sum\limits_{i=1}^{p}y^{(i)}\right\|-\left\|\sum\limits_{i=1}^{p}z^{(i)}-y^{(i)}\right\|\\&\geq \sum\limits_{i=1}^{p}\left\|y^{(i)}\right\|-\sum\limits_{i=1}^{p}\left\|z^{(i)}-y^{(i)}\right\|\\&\geq \sum\limits_{i=1}^{p}(\beta/2-2\epsilon 2^{-i})-\sum\limits_{i=1}^{p}2\epsilon 2^{-1}\geq \dfrac{p\beta}{2}-4\epsilon
\end{align*}
which is the desired inequality. Using this, it is easy to verify $(**).$ Indeed, taking $\xi_1$ and $\xi_2$ to be the points  $u_i$ and $v_i$ respectively in the intervals $[a_i,b_i ],$ and arbitrarily otherwise but the \textit{same} set of points in in the remaining intervals of $P,$ we have
\begin{align*}
\left\|(S(f,P,\xi_1)-S(f,P,\xi_2))\right\|=\left\|\sum\limits_{i=1}^{p}\dfrac{1}{N}z^{(i)}\right\|\geq \dfrac{1}{N}\left(\dfrac{p\beta}{2}-4\epsilon\right)\geq \dfrac{\alpha\beta}{4}.
\end{align*}
Finally, the sets $\{v_i;1\leq i\leq p\}\subset [a_i,b_i ]$ and $\{0=n_0\leq n_1\leq \cdots\}\subset\mathbb{N}$ having the afore-mentioned properties are constructed as follows. Since $\omega(f,u_2)\geq \beta,$ we can choose $v_1\in(a_1,b_1)$ such that $\|f(u_1 )-f(v_1 )\|\geq \beta/2.$ Also, there exists an integer $n_1>n_0$ such that $\sum_{j=n_1}^{\infty}\|z_{j}^{(1)}\|<\epsilon/2.$ Again since $\omega(f,u_2)\geq \beta,$ the continuity of  $e_{i}^{*}\circ f$ at $u_2$ for each $1\leq i\leq n_1$ yields $v_2\in (a_2,b_2)$ and $n_2>n_1$ such that
\begin{align*}
\sum\limits_{i=1}^{n_1}\|e_{i}^{*}\circ f(u_2)-e_{i}^{*}\circ f(v_2)\|<\epsilon/4\quad \text{and}\quad \sum_{j=n_2}^{\infty}\|z_{j}^{(2)}\|<\epsilon/4.
\end{align*}
Note that $e_{i}^{*}\circ f(u_2-v_2)=z_{i}^{(2)}$ and we get (a), (b) and (c) for $i=1, 2.$ Applying the same procedure to $u_i$ for $i=3, \cdots, p$ completes the construction. This completes the proof.

It turns out that $\ell_1$ crops up as an ubiquitous object in situations involving Riemann integration in Banach spaces. The first signs of this phenomenon are already evident in (b) above and together with the example in (c), the two can be subsumed within the class of the so called asymptotic $\ell_1$ Banach spaces defined below. However, some background material is needed to understand this definition.

Given a sequence $\left\{e_n \right\}_{n=1}^{\infty}$ in $X,$ a (non-zero) vector of the form $x=\sum_{i=m}^{n}a_ie_i,$ $\{a_i \}_{i=m}^{n}\subset \mathbb{R},$ shall be called a (\textit{block}) vector with respect to the given sequence  $\{e_n \}_{n=1}^{\infty}.$ The support of $x,$ denoted $\text{supp}(x),$ is defined to be the set of all integers $i$ for which $a_i\neq 0.$ Given two blocks $x$ and $y,$ we shall write $x<y$ if $\max(\text{supp}\,\, x)<\min(\text{supp}\,\, y)$  and that blocks $x_1,\cdots, x_n$ will be called successive if $x_1< . . . <  x_n.$ If $\{e_n \}_{n=1}^{\infty}$ is assumed to be a basic sequence, we shall say that a basic sequence $\{u_n \}_{n=1}^{\infty}$ of $X$ is a \textit{spreading model} of $\{e_n \}_{n=1}^{\infty}$ if there exist positive numbers $\epsilon_n\downarrow o$ such that for all $\{a_i \}_{i=1}^{\infty}\subset[-1,1],$ we have 
\begin{align*}
\left|\left\|\sum\limits_{i=1}^{n}a_ie_{k_{i}}\right\|-\left\|\sum\limits_{i=1}^{n}a_iu_{i}\right\|\right|<\epsilon_n,\quad\forall\,n\geq 1, 
\end{align*}
whenever $n\leq k_1< \cdots<  k_n.$
\begin{defn}
A Banach space $X$ is said to be an \textit{asymptotic} $\ell_1$ space with respect to a (normalised) basic sequence  $\{e_n \}_{n=1}^{\infty}$ if there exists $C > 1$ such that for each $n\in\mathbb{N},$ there exists a function $F_n:\mathbb{N}_{0}\rightarrow\mathbb{N}$ with $F_n (k)\geq k$ for all $k$ so that 
\begin{align*}
C^{-1}\sum\limits_{i=1}^{n}|a_i|\leq \left\|\sum\limits_{i=1}^{n}a_{i}x_{i}\right\|
\end{align*}
for all $\{a_i \}_{i=1}^{n}\subset\mathbb{R}$ and for all normalised successive blocks $\{x_i \}_{i=1}^{n}$ with respect to $\{e_n \}_{n=1}^{\infty}$ satisfying 
$$ F_n(0)\leq\min \text{supp}\,\, x_1\,\,\,\text{and}\,\,\,F_{n}(\max\text{supp}\,\, x_i)<\min\text{supp}(x_{i+1}),\quad i=1,2,\cdots,n-1. $$
\end{defn}
With this background, we are now in a position to state the following theorem of K. M. Naralenkov \cite{3} which reveals the role of $\ell_1$ in situations involving (RL)-property and at the same time, unifies the examples of Examples \ref{e1}. 

\begin{thm}[\cite{13}]
\begin{enumerate}
\item[(a):] Let $X$ be a Banach space having (RL)-property. Then each spreading model in $X$ is equivalent to the unit vector basis of $\ell_1.$
\item[(b):] Let $X$ be an asymptotic $\ell_1$ space with respect to a normalised basic sequence. Then $X$ has the (RL) property.
\end{enumerate}
\end{thm}
Besides of course the space $\ell_1$ and the Tsirelson space $T$ covered under Part (b) of the above theorem, other classes of Banach spaces falling under this category include the $\ell_1$-direct sum $(\sum_{n=1}^{\infty}E_n)_1$ of finite dimensional spaces $E_n$ as also the modified Tsirelson spaces $T_\theta (0< \theta <1).$ Note that the original Tsirelson space $T$ corresponds to $\theta=1/2:T=T_{1/2}.$ On the other hand, there are a whole lot of Banach spaces amongst the classical Banach spaces failing the (RL)-property. In fact, each of the spaces listed below fails the (RL)-property (See \cite{7}):
\begin{enumerate}
\item[(i)] $c_0,c,\ell_{\infty},C[a,b],L_1[a,b],L_{\infty}[a,b].$
\item[(ii)] $X^{*}$ if $X$ contains a copy of $\ell_1.$
\item[(iii)] Hilbert spaces. More generally,
\item[(iv)] Infinite dimensional Banach spaces which are uniformly convex.
\end{enumerate}
Thus we see that except for $\ell_1,$ no other space from amongst the classical Banach  spaces listed above figures in the list of spaces satisfying the (RL)-property. This motivates the case for weakening the definition of (RL)-property so as to include many more spaces within a larger class of Banach spaces defined by this weaker version of Riemann-Lebesgue property.
\begin{defn}[\cite{22}]
A Banach space $X$ is said to have \textit{Weak Riemann-Lebesgue (WRL) - property} if $f:[a,b]\rightarrow X$ is weakly continuous a.e on $[a, b]$ whenever it is Riemann integrable. 
\end{defn}
It can be shown that every Banach space with a separable dual has (WRL)-property (\cite{22}). However, $C[0,1]$ does not have the (WRL)-property as was shown by R.Gordon \cite{7}. The fact that separability of the dual is not a necessary condition for (WRL)-property was proved by C. Wang and K. Wan as the following result shows.
\begin{thm}[\cite{22}]
For a given measurable space $(\Omega,\sum)$ which is totally finite, complete and countably generated, the space $L_1(\Omega,\sum)$ has (WLP).
\end{thm}
As a far reaching generalisation of this result applied to the special case of $\Omega=[0,1]$ equipped with the Lebesgue measure, we have the following theorem of J. M. Calabuig et al :
\begin{thm}[\cite{5}]
Let $X$ be a Banach space having Radon-Nikodym property and let  $X^{*}$ be separable. Then the space $L_1([0,1],X)$ of Bochner integrable functions has (WLP). 
\end{thm}
\begin{cor}[\cite{5}]
Let $X$ be a separable reflexive Banach space. Then $L_1([0,1],X)$ has (WLP).
\end{cor}
\section{Certain Exotic integration phenomena in Quasi Banach spaces}\label{s2}
Let $X$ be a quasi Banach space, i.e., $X$ is a vector space equipped with a quasi norm  which is an $\mathbb{R}^{+}$- valued function on $X$ satisfying the properties of a norm, except that the triangle inequality comes with a constant greater than 1:
$$ \|x+y\|\leq k(\|x\|+\|y\|),\qquad\text{for all}\quad x,y\in X. $$
Amongst the familiar examples are the spaces $\ell_p$ and $L_p [a,b],\,\,\text{for}\,\,  0<p<1.$ The (linear) topology induced by a quasi norm is always non-locally convex, unless $k=1.$ Let us record below some exotic phenomena involving Riemann integration in quasi Banach spaces which are in sharp contrast to the analogous situation in the presence of local convexity. The contrast is especially striking in respect of the problem of primitives and of the fundamental theorem of calculus which, unlike in the case of Banach spaces as seen in Section \ref{s1}, break down in the quasi Banach space setting as will be seen presently. 

On the other hand, while three of the four fundamental theorems of functional analysis continue to hold good also in the quasi Banach space setting, it is chiefly the failure of Hahn Banach theorem in non-locally convex spaces that some of the familiar results from classical analysis extending even to the class of Banach space are no longer valid in the absence of local convexity. In fact, it turns out that in most of the cases, the validity of these results characterises local convexity of the space in question-the equivalence of Hahn Banach theorem with local convexity provides one such example. Let us also remark that the property of being quasinormed is essentially an infinite dimensional phenomenon since all Hausdorff linear topologies on a finite dimensional space are equivalent. To give a foretaste of what life would be like without local convexity, we provide below a sample of results which may be compared and contrasted with their corresponding counterparts in the familiar framework of Banach spaces. Most of these results have been proved in a series of important papers \cite{1}, \cite{2} and \cite{3}, besides many more by $F.$ Albiac and his collaborators. In the statements that follow, $X$ is assumed to be a (non-locally convex) infinite dimensional quasi Banach space.
\begin{thm}[Mazur-Orlicz \cite{12}]\label{t2.1}
Continuity of a function $f:[a,b]\rightarrow X$ does not imply Riemann integrability of $f.$ In fact, there always exists a continuous function $f:[a,b]\rightarrow X$ which fails to be Riemann integrable on $[a,b]$ $($which can even be chosen to be Riemann integrable on $[a,b-\epsilon]$ for each $\epsilon>0$ $)$.  
\end{thm}
\begin{thm}[Albiac and Ansorena \cite{1}]\label{t2.2}
 For $X$ such that $X^{*}$ is separating, there exists a continuous function $f:[a,b]\rightarrow X$ failing to have a primitive. $($In particular, this holds for$ X =l_p  , 0<p<1.)$ 
\end{thm}
An important ingredient in the proof of the above assertion is the following proposition which also shows that the mean value theorem does not hold for differentiable functions taking values in a quasi Banach space (with a separating dual). In fact, it can also be used to give an alternate proof of Theorem \ref{t2.1} (See also \cite[Theorem 3.5.1]{16}).
\begin{thm}[Albiac and Ansorena \cite{1}]\label{t2.3}
There always exists a continuously differentiable function $f$ on $[0, 1)$ into $X$ which does not admit a continuous extension to $[0, 1]$ but $f^{\prime}$ does!
\end{thm}
Thus, to produce a function $f$ with the indicated properties as asserted in Theorem \ref{t2.2}, let $f$ be the function as guaranteed by Theorem \ref{t2.3}, so that we can assume that $f^{\prime}$ is defined and continuous on $[0,1].$ Let $g= f^{\prime}$ on $[0,1]$ and assume that there exists a differentiable function $G:[0,1]\rightarrow X$ such that $G^{\prime} (t)=g(t)$ for all $t\in[0,1].$ Then $(f-G)^{\prime} (t)=0$ for all $t\in[0,1).$ By Theorem \ref{t2.6} quoted below, there exists $c\in\mathbb{R}$ such that $f(t) = G(t) + c$ for all $t\in[0,1),$ yielding that $f$ can be extended to a continuous function on $[0,1],$ thus contradicting the conclusion of Theorem \ref{t2.3}.
\begin{thm}[\cite{8}]\label{t2.4}
 For $X$ having a trivial dual, there exists a non-constant differentiable function $f:[a,b]\rightarrow X$ such that  $f^{\prime} (t)=0$ for all $t\in[a,b].$ 
\end{thm}
A special case of the above theorem is provided by the following simple example (due to Rolewicz) of a  function 
$f:[a,b]\rightarrow X$ for $X=L_p[a,b]$ $(0<p<1).$

Indeed, let $f$ be defined by $f(t)=\chi_{[0,t]}.$ We observe that
\begin{align*}
\left\|f(t+h)-f(t) \right\|=|h|^{1/p}
\end{align*}
and so 
\begin{align*}
\left\|f^{\prime}(t)\right\|=\left\|\frac{f(t+h)-f(t)}{h} \right\|=|h|^{\frac{1}{p}-1}.
\end{align*}
As $h\rightarrow 0$ and $p<1,$ the (RHS) goes to zero and so $f^{\prime}(t)=0$ for all $t\in[0,1].$

Amongst the positive results valid in the quasi Banach space setting, let us start with the following theorem due to Albiac and Ansoren \cite{2} which may be compared with Theorem \ref{t2.4} above
\begin{thm}[\cite{2}]\label{t2.6}
Assume that $X$ has a separating dual and let $J$ be an interval of real numbers such that $f:J\rightarrow X$ is differentiable with $f^{\prime} (t)=0$ for all $t\in J.$ Then $f$ is constant on  $J.$
\end{thm}
The following theorem was proved by M, Popov \cite{15}.
\begin{thm}\label{t2.5}
Let $f:[a,b]\rightarrow X$ be Riemann integrable. Then the function $F:[a,b]\rightarrow X.$ $($indefinite integral of $f)$ defined by 
$$ F(x)=\int\limits_{a}^{x}f(s)ds $$
is $($uniformly$)$ continuous on $[a,b].$ 
\end{thm}

\begin{thm}[Kalton \cite{9}]\label{t2.7}
For $X$ such that $X^{*}=(0),$ continuity of $f$ implies that $f$ has a primitive. $($In particular, this holds for $X = L_p [a,b],0<p<1).$
\end{thm}
\begin{thm}[Albiac and Ansorena \cite{3}]\label{t2.8}
Let $X$ be a (quasi) Banach space with a separating dual and let $f:[a,b]\rightarrow X$ be differentiable on $[a, b]$ so that $f^{\prime}$ is Riemann integrable on $[a, b].$ Then
$$ \int\limits_{a}^{b}f^{\prime}(t)dt=f(b)-f(a). $$
\end{thm}
In particular, such an f admits a primitive. Further, taking $X$ to be a Banach space yields (FTC) for Riemann integrable functions valued in a Banach space(See Section \ref{s1} (c)).

Let us close this section with a sample of proofs of some of these statements. We shall begin with the proof of Theorem \ref{t2.5}, which in the setting of Banach space is proved via the mean value property of Riemann integrals but which fails to hold in quasi Banach spaces. In the following proof due to Popov \cite{15}, we shall simplify the argument by assuming $k=1.$ This is possible because estimates involved in the proof can be appropriately scaled to absorb $k,$ so as to yield $k$-free estimates in the end.

Thus, to prove the continuity of $F$ on $[a,b],$ fix $t_1\in[a,b].$ For a given $ϵ>0,$ choose $\delta_1>0$ and $\delta_2>0$ such that the following estimates hold:
\begin{equation}\label{1.2}
\|\alpha f(t_1)\|<\epsilon/4,\qquad,\text{for all}\quad |\alpha|<\delta_1
\end{equation}
\begin{equation}\label{1.3}
\left\|S(f,P,\xi)-\int\limits_{a}^{b}f(s)ds\right\|<\frac{\epsilon}{8}
\end{equation}
for all partitions $P$ with $\|P\|<\delta_2.$ Here, $S(f,P,\xi)$ is the Riemann sum corresponding to the tagged partition $P=\{s_i,[t_{i-1},t_{i}],1\leq i\leq j \}$ of $[a,b]:$
$$ S(f,P,\eta)=\sum\limits_{i=1}^{j}f(s_i)\Delta_{i}\qquad(\Delta_i=t_i-t_{i-1}). $$
Choose $t\in[a,b]$ such that $|t-t_1 |<\delta=\min(\delta_1,\delta_2).$ Assume that $t_1<t$ (argument for the reverse inequality is analogous) and choose $\delta_3>0$ so that 
\begin{equation}\label{1.4}
\left\|S(f,P,\xi)-\int\limits_{t_{1}}^{t}f(s)ds\right\|<\dfrac{\epsilon}{2}
\end{equation}
for all partitions $P$ (of $[t_1,t]]$) with $\|P\|<\delta_3.$ Let $P_0 = \{s_i,[t_{i-1},t_i ],1\leq i\leq j\}$ be a partition of $[t_1,t]$ with $\|P_0\|<\delta_0=\min(\delta,\delta_3)$ and let $\xi$ denote the set $\{s_i \}_{i=1}^{j}.$  \\
\noindent Claim: $\|\sum_{i=1}^{j}f(s_i)\Delta_{i}\|<\epsilon/2+\|(t-t_1)f(t_1)\|.$ 

\noindent Using $P_0,$ we define tagged partitions $P_1$ and $P_2$ of $[a, b]$ as follows: 

\noindent$P_1$ is obtained from $P_0$ by adding to$\{t_{i}\}_{i=0}^{j}$ appropriately chosen points of $[a,b]\setminus [t_1,t]$ so that $P_1$ satisfies $\|P_1\|<\delta_0\leq \delta ,$ whereas  $P_2$ is obtained from $P_1$ by ignoring the points $t_1 ,\cdots,t_{j-1}.$ Since $|t-t_1 |<\delta,$ we see that  $\|P_2\|<\delta.$ To get tagged partitions, we choose the set $\xi_1$ of tags for $P_1$ as follows: For the intervals $[t_{i-1},t_i ],$ we retain the points $s_i$ as before and choose the left end points of the remaining intervals. In the case of $P_2,$ the set $\xi_2$ of tags is chosen to consist of the left end points of the intervals from $P_2.$ With this choice of tags, we have 
\begin{equation}\label{1.5}
S(f,P_1,\xi_1)-S(f,P_2,\xi_2)=\sum\limits_{i=1}^{j}f(s_i)\Delta_i-f(t_1)(t-t_1).
\end{equation}
On the other hand, since $\|P_1\|<\delta$ and $\|P_2\|<\delta,$ by \eqref{1.3} we have
\begin{align}\nonumber\label{1.6}
\|S(f,P_1,\xi_1)-S(f,P_2,\xi_2)\|&\leq \left\|S(f,P_1,\xi_1)-\int\limits_{a}^{b}f(s)ds\right\|\\\nonumber&\qquad+\left\|S(f,P_1,\xi_1)-\int\limits_{a}^{b}f(s)ds\right\|\\&\leq \frac{\epsilon}{8}+\frac{\epsilon}{8}=\frac{\epsilon}{4}.
\end{align}
Using \eqref{1.2}, \eqref{1.5} and \eqref{1.6}, we have
\begin{align*}
\left\|\sum\limits_{i=1}^{j}f(s_i)\Delta_i\right\|&\leq \left\|\sum\limits_{i=1}^{j}f(s_i)\Delta_i-f(t_1)(t-t_1)\right\|+\left\|f(t_1)(t-t_1)\right\|\\&<\frac{\epsilon}{4}+\frac{\epsilon}{4}=\frac{\epsilon}{2}
\end{align*}
Finally, combining the above estimate with \eqref{1.4}, we get
\begin{align*}
\|F(t)-F(t_1)\|&=\left\|\int\limits_{t_1}^{t}f(s)ds\right\|\\&\leq \left\|\int\limits_{t_1}^{t}f(s)ds-\sum\limits_{i=1}^{j}f(s_i)\Delta_i\right\|+\left\|\sum\limits_{i=1}^{j}f(s_i)\Delta_i\right\|\\&< \frac{\epsilon}{2}+\frac{\epsilon}{2}=\epsilon
\end{align*}
This completes the proof.

Let us now proceed to Theorems \ref{t2.6} and \ref{t2.8} which have been singled out in view of the similarity of argument underlying their proofs. In fact, we show Theorem \ref{t2.6} holds for an arbitrary interval $I$ in place of $[a, b].$ Thus, let us assume that $f(t)\neq f(s)$ for some $t, s \in I$ and let $f(t)\neq 0.$ Let $Y$ be the closed subspace of $X$ generated by $\{f(r):r\in I\}.$ Thus $Y\neq(0)$ and, by the given hypothesis, there exists a non-zero $x^{*}\in X^{*}$ such that $x^{*} (f(t))\neq x^{*} (f(s)).$ Since $f$ is differentiable, we see that $x^{*}\circ f$ is differentiable on  $I$ with derivative at each $r\in I$ given by
\begin{align*}
(x^{*}\circ f)^{\prime}(r)&=\underset{h\rightarrow 0}{\lim}\dfrac{x^{*}\circ f(r+h)-x^{*}\circ f(r)}{h}\\&=x^{*}\left(\underset{h\rightarrow 0}{\lim}\dfrac{f(r+h)-f(r)}{h}\right)\\&=x^{*}(f^{\prime}(r))=0.
\end{align*}
Finally, by the fundamental theorem of calculus applied to $(x^{*}\circ f)^{\prime}$ on $[s,t],$ we get
\begin{align*}
x^{*}\circ f(t)=x^{*}\circ f(s)+\int\limits_{s}^{t}x^{*}\circ f^{\prime}(r)dr=x^{*}\circ f(s),
\end{align*}
a contradiction.

A similar argument applies to prove Theorem \ref{t2.8}. Indeed, as noted above, for each $x^{*}\in X^{*},$ the composite map $x^{*}\circ f:[a,b]\rightarrow\mathbb{R}$ is differentiable with derivative $(x^{*}\circ f)^{\prime}(t)=x^{*}\circ f^{\prime}(t),\,t\in [a,b].$ Invoking (FTC) for real-valued functions(See Section \ref{s1}(b)), we have
\begin{align*}
x^{*}\left(\int\limits_{a}^{b}f^{\prime}(t)dt\right)&=\int\limits_{a}^{b}x^{*}\circ f^{\prime}(t)dt\\&=\int_{a}^{b}(x^{*}\circ f)^{\prime}(t)dt\\&=(x^{*}\circ f)(b)-(x^{*}\circ f)(a)\\&=x^{*}( f(b)- f(a)).
\end{align*}
Finally, using that $X^{*}$ separates points of $X$ completes the proof.
\section{Integration with respect to weaker forms of continuity:}\label{s3}
We have seen in Sections \ref{s1} and \ref{s2} that whereas in the quasi Banach space setting continuity is too mild a condition to imply Riemann integrability, for Banach space-valued functions continuity seems to be sufficiently strong to yield Riemann integrability. The question whether continuity relative to weaker topologies on the range space would still ensure Riemann integrability was addressed for the first time by Alexiewicz and Orlicz \cite{4} who produced an example of a $c_0$-valued weakly continuous function on $[0,1]$ which fails to be Riemann integrable(See also \cite{7}). However, a complete characterisation of Banach spaces for which weakly continuous functions are always Riemann integrable was given by V. M. Kadets \cite{10}.
\begin{thm}\label{t3.1}
For a Banach space $X,$ a function $f:[a,b]\rightarrow X$ being weakly continuous implies $f$ is Riemann-integrable if and only if $X$ is a Schur space $($i.e., weakly convergent sequences in $X$ are norm convergent$)$.
\end{thm}
A further strengthening of the above result was obtained by C. Wang and Z. Yang \cite{21} in the form of the following theorem. We shall need the following definition.
\begin{defn}\cite{21}
Given a locally convex topology $\tau$ (not necessarily compatible) on Banach spaces, we shall say that a Banach space $X$ has the $\tau$- Schur property (or $X$ is a $\tau$-Schur space) if $\tau$-convergence of a sequence in $X$ implies its norm convergence.
\end{defn}
Thus Schur property corresponds to $\tau$-Schur property by choosing $\tau$ to be the weak topology on $X.$
\begin{thm}[\cite{21}]
For a Banach space $X$ and a locally convex topology $\tau$ defined on it, the following statements are equivalent:
\begin{enumerate}
\item[(i)] 	Each function $f:[a,b]\rightarrow X$ which is $\tau$- weakly continuous is Riemann-integrable.
\item[(ii)] $X$ is a $\tau$- Schur space. 
\end{enumerate}
\end{thm}
Regarding Riemann integrability of function taking values inside the dual of a Banach space when equipped with weak$^{*}$-topology, Kadets \cite{10} notes that weak$^{*}$-continuity is too strong a condition to yield Riemann integrability in an infinite dimensional Banach space.
\begin{thm}[\cite{10}]\label{t3.4}
For a Banach space $X,$ each weak$^{*}$-continuous function $f:[a,b]\rightarrow X^{*}$ is Riemann-integrable if and only if $X$ is finite dimensional.
\end{thm}
In other words, the property involving \textit{Riemann integrability of weak$^{*}$-continuous functions} is a finite dimensional property in the following sense: Let us call this property (RW$^*$).
\begin{defn}[\cite{19}]
A property (P) of Banach spaces is said to be a \textit{finite dimensional $($FD$)$- property}  if it holds in each finite dimensional Banach space but fails in each infinite dimensional Banach space.
\end{defn}
\begin{ex}\label{e3.6}
\begin{enumerate}
\item[(i)] Heine-Borel Property ( $B_X$  is compact).
\item[(ii)] $X^{*}=X$ (algebraic dual of X).
\item[(iii)] Completeness of the weak-topology on $X.$\\
In view of the role played by (FD)-properties in the context of Riemann integrability, it will be useful to spend some time on the theme of finite dimensionality in an infinite dimensional context (See \cite{19} for a detailed treatment of this phenomenon). To this end, let us note that an important (FD)-property is provided by considering the so-called Hilbert-Schmidt property of a Banach space. Here, by a Hilbert Schmidt space we mean a Banach space $X$ with the property that a bounded linear map acting between Hilbert spaces and factoring over $X$ is a Hilbert-Schmidt map. Some well known examples of Banach spaces with this property are: $c_0,\ell_\infty,\ell_1,C(K),L_\infty(\Omega).$ As will be seen presently, the following (FD)-property will turn out to be of special importance in view of its 'universal' character in the sense to be made precise shortly. 
\item[(iv)] A Banach space which is simultaneously Hilbert and Hilbert Schmidt is finite dimensional (and conversely).
\end{enumerate}
\end{ex}
We shall now briefly describe certain useful properties of (FD)-properties and then proceed to investigate the (FD)-property encountered in Theorem \ref{t3.4} in the light of these properties.

\textit{Three important features of (FD)-properties: }

 There are three important features involving (FD)-properties but which manifest themselves only in an infinite dimensional context. These three features involving a given finite dimensional property (P) derive from: 
 \begin{enumerate}
 \item[(a)] Size of the set of objects failing (P).
 \item[(b)] Decomposition/Factorisation property of (P). 
 \item[(c)] Fr\'{e}chet space analogue of (P). 
 \end{enumerate}
\noindent(a)\textit{\underline{\textbf{ Size of the set}}} 

Given an (FD)- property (P), it turns out that for a given infinite- dimensional Banach space $X,$ the set of objects in $X$ failing (P) is usually 'very big': it could be \\
topologically big(dense)\\
algebraically big(contains an infinite-dimensional space)\\
big in the sense of category(non-meagre),\\ 
big in the sense of functional analysis (contains an infinite-dimensional closed subspace). \\

The following results proved by the author on the size of the set of objects failing the indicated (FD)-property show that these sets are indeed large.

\begin{ex}\label{e3.7}
\begin{enumerate}
\item[(i)](\cite{19})  $M(X)\setminus M_{bv}(X)$ together with the zero element contains an infinite dimensional space. This set is even known to be non-meagre.
\item[(ii)](\cite{14}) $M([0,1],X)\setminus B([0,1],X)$ together with the zero element contains an infinite-dimensional space.
\end{enumerate}
\end{ex}
Here, $M(X)$ and  $M_{bv} (X)$  stand respectively for the space of countably additive $X$-valued measures (respectively of bounded variation) whereas M and B denote, respectively, the classes of McShane and Bochner integrable functions on $[0, 1]$ taking values in $X.$ It may be noted that equality of sets displayed in each of these examples is an (FD)-property.

In respect of the (FD)-property (RW$^*$), in \cite{18} we had posed the problem regarding the size of the set determined by the failure of this property in an infinite dimensional Banach space. A positive answer to this question has been provided in a recent work by G M Cervantes (Murcia, Spain) which has been posted on the arXiv (Oct.29, 2015). Curiously the proof makes use of the 'fat' Cantor set as will be seen to be the case, in what follows, while dealing with the third feature of an (FD)-property involving the Fr\'{e}chet space analogue of the property (RW$^*$). 

\noindent \textbf{Theorem} (G M Cervantes): Given an infinite dimensional Banach space $X,$ the set of $w^*$-continuous functions  $f:[a,b]\rightarrow X^*$ which are not Riemann integrable together with the (identically) zero function contains a closed infinite dimensional space.

\noindent(b)\textit{\underline{\textbf{Decomposition of an (FD)-property}}} 

The (FD)-property given in Example \ref{e3.6} (iv) above comes across as a universal (FD)-property in the following sense:
Given an (FD)-property (P), it turns out that we can write (P) as a 'sum' of properties (Q) and (R): (P) = (Q) $\wedge$(R) in the sense that a Banach space $X$ verifies
\begin{enumerate}
\item[(Q)] iff $X$ is Hilbertian and 
\item[(R)] iff $X$ is Hilbert-Schmidt.
\end{enumerate}
In other words, an (FD)-property lends itself to a 'decomposition' (factorisation) as a 'sum' of properties (Q) and (R) which are characteristic of Hilbertisability and Hilbert-Schmidt property, respectively.
\begin{ex}
Let $\Pi_2$ denote the ideal of 2-summing maps and let $\Pi_{2}^{d}$ be the ideal of bounded linear maps acting between Banach spaces having a 2-summing adjoint. It turns out that the equality $\Pi_2(X,Y)=\Pi_{2}^{d}(X,Y)$ for all Banach spaces $Y$ is an (FD)-property. A decomposition of this (FD)-property that suggests itself is the natural one involving the decomposition of the above equality into inclusion of each set appearing above inside the other:
\begin{enumerate}
\item[(a)] $\Pi_{2}(X,Y)\subset \Pi_{2}^{d}(X,Y),$ for all $Y.$
\item[(b)] $\Pi_{2}^{d}(X,Y)\subset \Pi_{2}(X,Y),$ for all $Y.$
\end{enumerate}
\end{ex}
What is indeed remarkable is that for a Banach space $X,$ (Q) holds if and only if $X$ is Hilbertian whereas (R) is valid precisely when $X$ is a Hilbert Schmidt space (See \cite{19} for details). Moreover, it suffices to choose $Y = \ell_2$ as a test space in (Q).

Regarding the (FD)-property (RW$^*$) under discussion, a possible decomposition of this property would entail the identification of a locally convex topology $\tau$ on the dual of a Banach space $X$ stronger than the weak$^{*}$ topology such that 

(Q) Each $f:[a,b]\rightarrow X^{*}$ continuous w r t the topology $\tau$ is Riemann integrable if and only if $X$ is a Hilbert space.\\
(R) Each weak$^{*}$-continuous map $f:[a,b]\rightarrow X^{*}$ is $\tau$- continuous if and only if $X$ is a Hilbert Schmidt space.

The above discussion suggests the following problem:
\begin{pr}\label{}
Describe the existence of a locally convex topology $\tau$ on the dual of a Banach space $X$ such the conditions (Q) and (R) as specified above in the last paragraph are satisfied.
\end{pr}
\noindent(c)\textit{\underline{\textbf{ Fr\'{e}chet space setting}}} 

When suitably formulated in the setting of Fr\'{e}chet spaces $X,$ it turns out that in most of the cases, there exist infinite dimensional Fr\'{e}chet spaces in which it is possible to salvage a given (FD)- property. It also turns out that, at least in most cases of interest, the class of Fr\'{e}chet spaces in which this holds coincides with the class of nuclear spaces (in the sense of Grothendieck) or, in some cases, with a class of spaces which are close relatives of nuclear spaces. The fact that a Banach space can never be nuclear unless it is finite dimensional shows that the nuclear analogue of an (FD)-property provides a maximal infinite dimensional setting in which the given (FD)-property can be saved. An important example to illustrate this situation involves the Fr\'{e}chet space analogue of the Dvoretzky-Rogers property which holds exactly when the underlying Fr\'{e}chet space is nuclear. Let us recall that the Dvoretzky-Rogers theorem is the statement that the Banach space analogue of the Riemann rearrangement theorem from classical analysis holds exactly when the Banach space in question is finite dimensional. The rest of this section is devoted to a brief sketch of the proof of the Fr\'{e}chet analogue of Theorem \ref{t3.4}, i.e., of the property (RW$^*$). To this end, we shall begin with the definition of Riemann integrability in Fr\'{e}chet spaces which is an obvious generalisation of the definition of Riemann integrability for Banach space -valued functions.
\begin{defn}
Let $X$ be a Fr\'{e}chet space and let 〖$\{p_m\}_{m=1}^{\infty}$ be a sequence of seminiorms generating the (Fr\'{e}chet)-topology of $X.$ We shall say that a function $f:[a,b]\rightarrow X$ is Riemann-integrable  if the following holds:\\
(*) $\exists\, x\in X$ such that $\forall\,\epsilon>0$ and $n\geq 1,$ $\exists \delta=\delta(\epsilon,n)>0$ such that for each tagged partition $P=\{s_i, [t_{i-1},t_i], 1\leq i\leq j \}$ of $[a,b]$ with
$$ \|P\|=\underset{1\leq i\leq j}{\max}(t_{i}-t_{i-1})<\delta, $$
we have
$$ p_{n} (S(f,P)-x)<\epsilon, $$
where, $S(f,P)$ is the Riemann sum of $f$ corresponding to the tagged partition $P=\{s_i, [t_{i-1},t_i]; 1\leq i\leq j \}$ of $[a,b]$ with $a=t_0<t_1<\cdots<t_j=b$  and $s_i\in [t_{i-1},t_i], 1\leq i\leq j.$ Here, the (unique) vector $x,$ to be denoted by 
$$ \int\limits_{a}^{b}f(t)dt, $$
shall be called the Riemann-integral of $f$ over $[a, b].$
\end{defn}
The following statement, proved by the author is the Fr\'{e}chet analogue of Kadet’s theorem stated above.
\begin{thm}\cite{18}\label{t3.11}
For a Fr\'{e}chet space $X,$ each $X^{*}$-valued weakly$^{*}$-continuous function is Riemann integrable if and only if $X$ is a Montel space.
\end{thm}
(A metrisable locally convex space is said to be a Montel space if closed and bounded subsets of X are compact).

Since Banach spaces which are Montel are precisely those which are finite dimensional, Theorem 3.11 yields Kadet’s theorem as a very special case.\\

\noindent\textit{\underline{\textbf{  Ingredients of the proof:}}} \\
 \textit{Construction of a 'fat' Cantor set.}\\
A 'fat' Cantor set is constructed in a manner analogous to the construction of the conventional Cantor set, except that the middle third subinterval to be knocked out at each stage of the construction shall be chosen to be of a suitable length $\alpha$ so that the resulting Cantor set shall have nonzero measure. In the instant case, each of the $2^{k-1}$ subintervals $A_k^{(i)}$  $( i=1,2,…,2^(k-1))$  to be removed at the kth stage of the construction from each of the remaining subintervals $B_{k}^{(i)}$ $( i=1,2,…,2^(k-1))$ at the (k-1)th stage shall be of length $\alpha=d(A_{k}^{(i))}=\frac{1}{2^{k-1}}\dfrac{1}{3^k},$ in which case $d\left(B_k^{(i)}\right)=\frac{1}{2^{k}}\left(1-\sum_{j=1}^{k}\frac{1}{3^j}\right)$ and, therefore, $d(C)=\frac{1}{2}.$ \\

\noindent\textit{Fr\'{e}chet analogue of Josefson-Nissenzwieg theorem:}

\noindent It is a well known theorem of Josefson and Nissenweig that from the unit sphere in the dual of an infinite dimensional Banach space, it is always possible to extract a sequence which is weak$^{*}$-null. A Fr\'{e}chet analogue of this important theorem and useful for our purpose was proved by Bonet, Lindsrtom and Valdivia (Two theorems of Josefson-Nissenweig type for Fr\'{e}chet spaces, Proc. Amer. Math. Soc. 117 (1993), 363-364):

\noindent\textbf{Theorem:} A Fr\'{e}chet space $X$ is Montel if and only if each weak$^{*}$-null sequences in $X^{*} $ is strong$^{*}$-null.

\noindent\textit{Sketch of proof of Theorem \ref{t3.11}}:\\
\textit{Necessity}: This is a straightforward consequence of Theorem (b) above combined with the sequential completeness of $X^{*}$ in its weak$^{*}$ topology.\\
\textit{Sufficiency}: Assume that $X$ is not Fr\'{e}chet Montel. By (b), there exists a sequence in $X^{*}$ which is weak$^{*}$-null but not strong$^{*}$-null. Denote this sequence by $\{x_{n}^{*} \}_{n=1}^{\infty}.$ Write $A_{k}^{(i)}=[a_{k}^{(i)},b_{k}^{(i)}]$ and choose a function $\phi_{k}^{(i)}:[0,1]\rightarrow\mathbb{R}$ which is piecewise linear on $A_{k}^{(i)}$ and vanishes off $A_{k}^{(i)}$  with $\|\phi_{k}^{(i)}\|\leq 1.$
Put
\begin{align*}
h_{k}(t)=\sum\limits_{i=1}^{2^{k-1}}\phi_{k}^{(i)}(t),\,t\in [0,1],
\end{align*}
and define
\begin{align*}
f(t)=\sum\limits_{i=1}^{2^{k-1}}h_{k}(t)x_{n}^{*},\,t\in [0,1].
\end{align*}
Claim 1: $f$ is weak$^{*}$-continuous.

Note that each $h_k$ is continuous on $[0,1]$ such that $\|h_k\|\leq 1.$ The proof of the claim is achieved by showing that the series defining f is uniformly convergent in $X_{\sigma}^{*}.$ Indeed, fix $\epsilon>0$ and $x\in X.$ Choose $K_0$ such that $|\langle x_{k}^{*},x\rangle|<\epsilon$ for all $k\geq K_0.$ By the definition of $h_k,$ it follows that
\begin{align*}
\sum\limits_{n=k+1}^{\infty}h_{n}(t)x_{n}^{*}=0\quad\text{for}\quad t\in \bigcup_{n=1}^{k}\bigcup_{n=1}^{2^{n-1}}A_{n}6{(i)}
\end{align*}
and that for each $t\in\bigcup_{n=k+1}^{\infty}\bigcup_{i=1}^{2^{n-1}}A_{k}^{(i)},$ we can choose $k_0>k$ such that
\begin{align*}
\sum\limits_{n=k+1}^{\infty}h_{n}(t)x_{n}^{*}=h_{k_{0}}(t)x_{k_{0}}^{*}. 
\end{align*}
It follows that for all $k>K_{0}$ and for all $t\in [0,1],$ we have 
\begin{align*}
\left|\langle f(t)-\sum\limits_{n=1}^{k}h_{n}(t)x_{n}^{*},x\rangle\right|&=\left|\langle\sum\limits_{n=k+1}^{\infty}h_{n}(t)x_{n}^{*},x\rangle\right|\\&\leq \left|\langle h_{k_{0}}(t)x_{k_{0}}^{*},x\rangle\right|\leq \left|\langle x_{k}^{*},x\rangle\right|<\epsilon
\end{align*}
Claim 2: f is not Riemann integrable.\\
Here we use the fact that the Cantor set $C$ constructed above has measure equal to $1/2$ and then produce a bounded subset $B$ of $X$ and for each $\delta>0$ tagged partitions $P_1$ and $P_2$ of $[0, 1]$ with $\|P_1\|<\delta,$ $\|P_2\|<\delta$ such that 
\begin{align*}
p_{B}(S(f,P_1)-S(f,P_2))>1/2,
\end{align*}
where  $p_B$ is the strong$^{*}$-seminorm on $X^{*}$  corresponding to $B$ defined by
\begin{align*}
p_{B}(f)=\underset{x\in B}{\sup}|f(x)|,\quad f\in X^{*}.
\end{align*}
The desired bounded set $B$ is obtained by the assumed hypothesis that the sequence $\{x_n^{*} \}_{n=1}^{\infty}$ is not strong$^{*}$-null which means that, passing to a subsequence if necessary,  $p_B (x_{n}^{*} )>1$ for some bounded subset $B\subset X$. Fix $\delta>0$ and choose $m\geq 1$ such that $2^{-m}<\delta.$ Note that $d(B_{m}^{(i)})<2^{-m}$ for $m \geq 1$ and for $1\leq i\leq 2^{m-1}.$ The desired partitions $P_1=\{(s_{i},[t_{i-1},t_i ]); 1\leq i\leq N_m\}$ and $P_2=\{(s_{i}^{\prime},[t_{i-1},t_i ]); 1\leq i\leq N_m\}$  with $\|P_1\|<\delta,$ $\|P_{2}\|<\delta$ consisting of the same number $N_m$  of intervals are obtained by demanding that
\begin{enumerate}
\item[(a)] Both partitions contain the sets $B_{m}^{(i)}$ for $i =1,\cdots,2^{m-1}.$
\item[(b)] $t_{i}-t_{i-1}<2^{-(m-1)}$ for $1\leq i\leq N_{m}.$
\item[(c)] $s_{i}=s_{i}^{\prime}$ if $[t_{i-1},t_{i}]\neq B_{m}^{(i)}, i=1,\cdots,2^{m-1}.$
\item[(d)] $s_i=c_{m}^{i}$ and $s_{i}^{\prime}=a_{m}^{i},$ if $[t_{i-1},t_i]=B_{m}^{(i)},$ $i=1,\cdots,2^{m-1},$ where $c_{m}^{i}$ is chosen to be the midpoint of $A_{m}^{(i)}$ such that $\phi_{m}^{(i)}(c_{m}^{(i)})=1.$
\end{enumerate}

Now (c) gives $f(s_i)=f(s_{i}^{\prime})$ if $[t_{i-1},t_i ]\neq B_{m}^{(i)}$ whereas (d) yields
$$ f(s_i)=h_{m}(s_i)x_{m}^{*}=\phi_{m}^{(i)}(c_{m}^{(i)})x_{m}^{*}=x_{m}^{*},\quad f(s_{i}^{\prime})=0,\quad\text{if}\quad[t_{i-1},t_i]=B_{m}^{(i)}. $$
Finally, the above construction gives
\begin{align*}
p_{B}(S(f,P_1)-S(f,P_2))&=p_{B}\left(\sum\limits_{i=1}^{N_m}(f(s_i)-f(s_{i}^{\prime}))\right)(t_i-t_{i-1})\\&=p_{B}\left(\sum\limits_{i=1}^{2^{m-1}}(h_{m}(s_i)x_{m}^{*}d(B_{m}^{(i)}))\right)\\&=p_{B}(x_{m}^{*})\sum\limits_{i=1}^{2^{m-1}}d(B_{m}^{(i)}))\\&=p_{B}(x_{m}^{*})2^{m-1}\left[2^{-(m-1)}\left(1-\sum\limits_{j=1}^{m}\dfrac{1}{3^{j}}\right)\right]>\frac{1}{2}.
\end{align*}
We conclude with the following problems which appear to be open.
\begin{pr}
Characterise the class of Banach spaces $X$ such that weakly-continuous functions $f:[a,b]\rightarrow X$ have a primitive $F:$
$$ F^{\prime}(t)=f(t),\qquad\forall\,t\in [a,b], $$
i.e.,
\begin{align*}
\underset{h\rightarrow 0}{\lim}\left\|\dfrac{F(t+h)-F(t)}{h}-f(t)\right\|=0,\quad \forall t\in [a,b].
\end{align*}
\end{pr}
It is not clear if the class of spaces enjoying this property contains an infinite dimensional Banach space. However, if the existence of an almost primitive (entailing differentiability of F a.e.) is demanded instead of a primitive, the obvious examples meeting this requirement are provided by the class of Banach spaces satisfying the Schur property and the (RL)-property ($\ell_1,$ for example). This motivates the following problem.
\begin{pr}\label{p3.13}
Describe the class of Banach spaces $X$ such that each Riemann integrable function $f:[a,b]\rightarrow X$ is differentiable a.e.
\end{pr}
Here it may be useful to recall that whereas for Bochner integrable functions the indefinite integral is differentiable almost everywhere, in the case of infinite dimensional Banach space $X,$ there always exists an $X$-valued Pettis (even McShane) integrable function whose indefinite integral is non-differentiable on a set of positive measure (see \cite{11})! The fact that there exist infinite dimensional Banach spaces having (RL)-property shows that in the foregoing statement, it is not possible to choose the function to be Riemann integrable.

 We conclude our discussion of primitives with the 'a.e.-analogue' of the fundamental theorem of calculus by asking if it is possible to recover f from its almost primitive $F$ (as a definite integral). The answer is provided by the following theorem of C. Volintiro \cite{20}:
\begin{thm}
Suppose that for a Banach space $X$ and a (Borel) measurable Riemann integrable function $f:[a,b]\rightarrow X,$ there exists a continuous almost primitive $F$ of $f,$ i.e., a continuous function $F:[a,b]\rightarrow X$ such that $F^{\prime}(t)$ exists and $F^{\prime} (t)=f(t)$ outside a Lebesgue null set $B\subset[a,b].$  Further assume that $F(B)$ has Hausdorff measure equal to zero. Then
\begin{align*}
F(x)=\int\limits_{a}^{x}f(t)dt,\quad x\in [a,b].
\end{align*}
\end{thm}
\begin{flushleft}
\textbf{Epilogue}
\end{flushleft}
\noindent Some Related Developments:

 \noindent We have tried to provide a smorgasbord of results pertaining to vector-valued Riemann integration, with the functions involved taking values in various classes of spaces: Banach spaces, quasi Banach spaces and Fr\'{e}chet spaces, with the main aim to draw attention to the rich interplay of ideas between certain important aspects of vector-valued Riemann integration theory and the geometry/structure of the spaces in question. In the context of these classes of spaces, it is mainly within the framework of Banach spaces that the theory involving vector integration has been pursued vigorously and successfully, marked by certain spectacular developments that have been witnessed in recent years. We conclude with a brief overview (without technical details)  of these developments. 
\begin{flushleft}
(i) A litany of integration theories
\end{flushleft}
In the previous paragraphs, we have made no attempt to discuss a whole lot of other vector-valued integration theories that are in vogue right now, nor the rich theory underlying the question involving the description of Banach spaces $X$ witnessing coincidence of various integration theories of functions valued in $X.$ Apart from vector-valued Riemann integration that has been the main thrust of this paper, other integration theories having been proposed in the setting of Banach spaces invariably reduce to Riemann or Lebesgue integration in the finite dimensional context. Amongst the more prominent ones of these theories-besides, of course, vector-valued Riemann integration (R) - are the theories known after Bochner (B), Riemann-Lebesgue (RL), Birkhoff (BR), McShane (M), Henstock-Kurzweil (HK), Pettis (P). The relationship between these integration theories is provided by the following implication diagram:
 
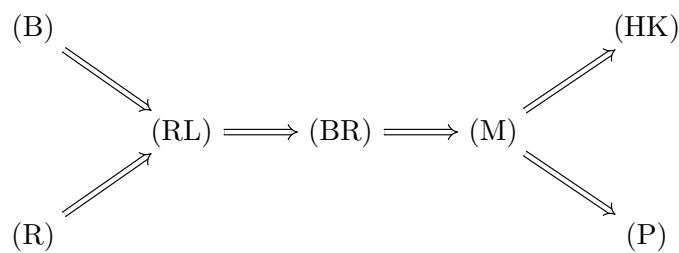
\[
\begin{tikzcd}[%
    ,every cell/.append style={align=center}
    ,every arrow/.append style={Rightarrow}
    ]
\text{(B)}\drar & & & & \text{(HK)}\\
& \text{(RL)}\rar & \text{(BR)}\rar & \text{(M)}\urar\drar & \\
\text{(R)}\urar & & & & \text{(P)}
\end{tikzcd}
\] 
Further, there are examples to show that the above implications are strict. This motivates the problem of describing necessary and sufficient conditions for the reverse implications to hold. Unfortunately, not much is known in this direction except in a few cases including the implication  $\text{(P)}\Rightarrow \text{(M)}$ which is known to be valid for Banach spaces that are weakly compactly generated (WCG)(A. Aviles et al, Jour. Func. Anal. 259(11), 2010, 2776-2792). Again, it is folklore that the equivalence $\text{(P)}\Rightarrow \text{(B)}$ characterises finite dimensionality of the underlying Banach space as does the weaker property $\text{(P)}\Leftrightarrow\text{(M)}.$ Furthermore, a well-known theorem of D. H. Fremlin (Illinois. Jour. Math. 38(3), 1994, 471-479) tells us that $\text{(M)}\Leftrightarrow (\text{HK}+\text{P}).$ The remaining implications though, remain largely unexplored and so provide a fertile area for further research.  

An interesting insight into the problem of equivalence of various integration theories is provided by considering classes of operators that improve integrability. The first result along these lines was provided way back in 1972 when J. Diestel (Math. Ann. 196 (1972), 101-105) was able to characterise absolutely summing maps as those which (continuously) map (the space of) Pettis integrable functions into Bochner integrable ones. As a refinement of this result, J. Rodriguez (Jour. Math. Anal. Appl. 316(2), 2006, 579-600) showed that the same characterisation holds for McShane integrable functions in place of the class of Pettis integrable functions. Very recently, these ideas have been pursued in a series of papers by E. A. S.  Sanchez and his co-workers to the setting of the so-called (RS)-abstract $p$-summing maps.      

\begin{flushleft}
(ii) Lebesgue Differentiation Property(LDP) 
\end{flushleft} 
The classical Lebesgue differentiation theorem states that the indefinite integral of a Lebesgue integrable function is differentiable a.e. and that its derivative equals the given function a.e. The corresponding statement for a Riemann integrable function is a simple consequence of a.e. continuity of such a function. In case of vector integration - as pointed out in the discussion following Problem 3013 - the said property is known to hold for Bochner integrable functions, whereas there exist examples of Pettis (even McShane) integrable functions taking values inside 
$C[0,1]$ or $\ell_2$ and failing the (LPD). Further, in confirmation of a conjecture of G.E.F Thomas, V.M. Kadets \cite{11}  in 1994 showed that the indicated property fails in each infinite dimensional Banach space. This was followed by an important work of S.J. Dilworth and M. Girardi(Quaes. Math. 18(1995), 365-380) who provided a far reaching strengthening of Kadets result to yield for each infinite dimensional Banach space $X$ the existence of an $X$-valued Pettis integrable functions with an everywhere non-differentiable indefinite integral. In the process they also show how the cotype of the underlying Banach space is closely related to the degree of non-differentiability of the indefinite integral. In a recent work by B. Bongiorno, U.B.Darji and L. Di Piazza (Monat. Math. 177(3), 2015, 345-362) and in conformity with the philosophy as spelt out in Section $3$(a), the authors show that for each infinite dimensional Banach space $X,$ the set of $X$-valued Pettis (McShane) integrable functions on $[0,1]$ with primitives failing to be differentiable anywhere on $[0,1]$ even contains an infinite dimensional space! Again, the problem involving the validity/failure of the (LDP) in respect of integration theories other than those of Bochner, McShane and Pettis remains wide open. 
\begin{flushleft}
(iii) Limit Sums of non-integrable functions
\end{flushleft}
As opposed to what has been discussed in the previous paragraphs regarding integration of vector-valued functions, another fascinating area of research involving the issue of non-integrability pertains to the structure of 'limit sums' of non-Riemann integrable functions and their relationship with the structure of the underlying Banach spaces - a theme which we shall exclude from our discussion, considering that it has been dwelt upon at length in their beautiful monograph by M.I. Kadets and V.M. Kadets (Series in Banach spaces, Birkhauser Verlag, 1997). We conclude by quoting an important fact proved by V.M. Kadets as a sample result belonging in this circle of ideas which asserts that the convexity of the set of 'limit sums' in a Banach space characterises the $K$-convexity of the space in question.

% ------------------------------------------------------------------------
\end{document}